\documentclass[12pt]{amsart}
\usepackage{amssymb}
\usepackage{amsthm}
\usepackage{enumerate}
\usepackage[cp1251]{inputenc}
\usepackage[russian]{babel}

\renewcommand{\phi}{\varphi}

\newcommand{\eps}{\varepsilon}
\newcommand{\Th}{\Theta}
\newcommand{\Om}{\Omega}
\newcommand{\gS}{\mathfrak {S}}
\newcommand{\Kth}{K_\theta}
\renewcommand{\C}{\mathbb C}
\newcommand{\T}{\mathbb T}
\newcommand{\R}{\mathbb R}

\newcommand{\Z}{\mathbb Z}

\newcommand{\tU}{\tilde U}
\newcommand{\tS}{\tilde S}
\newcommand{\UU}{V}

\newcommand{\PW}{\mathcal{PW}}

\newtheorem{Thm}{Theorem}[section]
\newtheorem{Lem}[Thm]{Lemma}
\newtheorem{Prop}[Thm]{Proposition}
\newtheorem{Cor}[Thm]{Corollary}
\newtheorem{Rem}[Thm]{Remark}
\newtheorem{Ex}[Thm]{Example}

\def\beginpf{\noindent {\bf Proof.} \ }

\title[On perturbations of the semigroup of shifts]
{On perturbations of the isometric semigroup \\  
of shifts on the semiaxis}
\author{G.G.~Amosov, A.D.~Baranov, V.V.~Kapustin}
\date{}
\thanks{The work is partially supported by Federal program 2.1.1/1662,
by RFBR grant 08-01-00723, and by the President of Russian Federation 
grant NSH 2409.2008.1.}

\begin{document}
\maketitle
\sloppy
\noindent
{\small {\bf Abstract.} 
We study perturbations $(\tilde\tau_t)_{t\ge 0}$ of the semigroup of 
shifts $(\tau_t)_{t\ge 0}$ on $L^2(\R_+)$ with the property that 
$\tilde\tau_t - \tau_t$ belongs to a certain Schatten--von Neumann class 
$\gS_p$ with $p\ge 1$. We show that, for the unitary component 
in the Wold--Kolmogorov decomposition of the cogenerator of the semigroup 
$(\tilde\tau_t)_{t\ge 0}$, {\it any singular}
spectral type may be achieved by $\gS_1$ perturbations. 
We provide an explicit construction 
for a perturbation with a given spectral type
based on the theory of model spaces 
of the Hardy space $H^2$.
Also we show that we may obtain
{\it any} prescribed spectral type for the 
unitary component of the perturbed semigroup 
by a perturbation from the class $\gS_p$ 
with $p>1$.
 }

\bigskip

\section{Introduction}\label{vved}

Consider the isometric semigroup $(\tau_t)_{t\ge 0}$ of shifts 
on the space $L^2(\R_+)$,
$$
(\tau_t f)(x)= 
\begin{cases}
f(x-t), & x\ge t,\\
0,      & x<t,   \\
\end{cases}
\qquad f\in L^2(\R_+).
$$
In this paper we are concerned with perturbations
$(\tilde\tau_t)$ of the semigroup $(\tau_t)$ satisfying the 
following properties:

$(\tilde\tau_t)_{t\ge 0}$ is a strongly continuous semigroup 
of isometric operators on $L^2(\R_+)$;

the difference $\tilde\tau_t-\tau_t$ belongs to a certain
Schatten--von Neumann ideal $\gS_p$ for every $t>0$.

The central problem considered in this paper is to describe all possible
spectral types of perturbed isometric semigroups. 
The spectral type of a semigroup determines the semigroup uniquely up
to the unitary equivalence; it is defined by the spectral type
of the cogenerator of the group (see definition in \S2). 
For $p=1$ it follows from the stability 
of the absolutely continuous spectrum 
of the unitary dilation that 
the absolutely continuous parts of the cogenerators of the unitary dilations 
of the semigroups $(\tau_t)$ and  $(\tilde\tau_t)$ are unitarily equivalent 
(see \S2 for details). The cogenerator of the semigroup $(\tau_t)$ 
is unitarily equivalent to the unilateral shift operator 
on the Hardy space $H^2$. Thus, for $p=1$,
our problem reduces to the description 
of all possible singular parts, and we show that   
{\it any singular} type  may be realized by  some 
semigroup $(\tilde\tau_t)_{t\ge 0}$.
For $p=2$ it was shown in  \cite{AB2, AB3} 
that {\it any} spectral type of the unitary component is possible;
for the case of a singular spectral measure a model of such perturbations 
was constructed. Here we show that analogous results 
are valid for all $p>1$. 

A possible motivation for our study is a connection with Markovian
perturbations of the unitary group of shifts.
The semigroup $(\tau_t)_{t\ge 0}$ 
is the restriction to $L^2(\R_+)$ of the unitary group 
of shifts $(\gamma_t)_{t\in \R}$, $(\gamma_t f) (x) = f(x-t)$,
on the space $L^2(\R)$ on the whole line
(we identify the spaces $L^2(\R_+)$ and $L^2(\R_-)$ with the subspaces
of functions from $L^2(\R)$ identically equal to zero on the semiaxes $\R_-$
and $\R_+$, respectively). Consider a perturbed unitary group 
$(\tilde\gamma_t)_{t\in \R}$ and assume that it has
the so-called Markovian property, which means
that the perturbed operators coincide with the unperturbed ones on the 
left semi-axis $\R_-$ as $t<0$, i.e., 
\begin{equation}
\label{100}
\gamma_t f=\tilde\gamma_t f \quad\text{ if }\quad f=0 \;\text{ on }\; \R_+,
\qquad t<0.
\end{equation}
As usual, the Markovian property may be interpreted in the sense that
``the past does not depend on the future''. Markovian 
perturbations with the additional property 
$\tilde\gamma_t-\gamma_t\in \gS_2$, $t\in \R$, were investigated by the
first author (see, e.g., \cite{Amo00}) in connection with cocycle 
perturbations of the flow of Powers shifts \cite{Arv}.
\medskip

Sometimes it will be convenient to work with the model in 
the Hardy space $H^2$ on the unit circle $\T$ which is 
unitarily equivalent to the original one (see \S4).
The semigroup $(\tau_t)_{t\ge 0}$ on $L^2(\R_+)$ is unitarily equivalent
to the semigroup of operators on $H^2$ of multiplication by the functions
$\phi_t$, $t\ge 0$, where
\begin{equation}
\label{80}
\phi_t(z)=\exp\bigg(t\frac{z+1}{z-1}\bigg).
\end{equation}
In this case the cogenerator of the unperturbed semigroup is the unilateral shift $S$,
i.e., the operator of multiplication by $z$ in $H^2$.
Denote by $\tS$ the cogenerator of the perturbed semigroup; then its
elements are of the form $\phi_t(\tS)$.

Now we introduce a new parameter
of our model of perturbations, 
namely, an inner function $\theta$ in the unit disk
(a function $\theta\in H^2$ is said to be inner if $|\theta|=1$ almost
everywhere on $\T$). Consider the  $S$-coinvariant subspace 
$$
\Kth=H^2\ominus\theta H^2.
$$ 
For the theory the backward shift invariant subspaces,
also called model subspaces, see \cite{nk, Nik}.

In what follows we are interested in perturbations $\tS$ of
the shift operator $S$ with the following properties:

(i)\; $\tS$ is an operator on $H^2$, diagonal 
with respect to the decomposition $H^2=\Kth\oplus\theta H^2$;

(ii)\; $\tS$ acts on $\theta H^2$ as multiplication by $z$;

(iii)\; the restriction of $\tS$ to $\Kth$ is a unitary operator, 
and 1 is not its eigenvalue. 

Conditions (i)--(iii) mean that $\tS$ is an isometry, whose unitary
and completely non-unitary parts act on $\Kth$ and $\theta H^2$, respectively.
The number 1 does not belong to the point spectrum of $\tS$, 
and the isometric semigroup with cogenerator $\tS$ is well defined.

Now we state our main results. First of them shows that any singular 
summand $\UU$ may be obtained by means of 
trace class perturbations of the isometric 
semigroup of shifts in the model in the unit circle. 

\begin{Thm} \label{thm1}
Let $\UU$ be a singular unitary operator of multiplicity $n\leq\infty$
such that $1$ is not an eigenvalue of $V$, and let $\eps>0$. Then there exist
an inner function $\theta$ and an operator $\tS$ with the properties
{\rm (i)--(iii)} such that

$a)$ \, the restriction of $\tS$ to $\Kth$ is unitarily equivalent to $\UU$;

$b)$ \, $rank\,(\tS-S)\leq n$;

$c)$ \, $||\tS-S||_{\gS_1}\leq\eps$;

$d)$ \, $\phi_t(\tS)-\phi_t(S)\in\gS_1$ for all $t>0$.
\end{Thm}

As an immediate consequence of Theorem \ref{thm1}, we obtain
the following statement about perturbations of the semigroup of shifts 
on the semiaxis.

\begin{Thm} \label{thm8}
Let $(\tau_t)_{t\ge 0}$ be the semigroup of shifts on $L^2(\R_+)$.

$1)$ If $(\tilde\tau_t)$ is an isometric semigroup such that 
$\tilde\tau_t-\tau_t\in\gS_1$ for all $t\ge 0$, then the 
cogenerator of the semigroup $(\tilde\tau_t)$ 
is unitarily equivalent to the direct sum 
of the shift operator (of multiplicity one) 
and a unitary operator with a singular spectral measure.

$2)$ For any singular unitary operator $V$ with the only restriction that 
$1$ does not belong to its point spectrum, there exists 
an isometric semigroup $(\tilde\tau_t)$ whose cogenerator 
is unitarily equivalent to the operator $S\oplus V $, and 
$\tilde\tau_t-\tau_t\in\gS_1$ for all $t\ge 0$.
\end{Thm}

As another corollary of Theorem \ref{thm1} we show that 
for $p>1$ an analogous statement is true without the assumption
that the unitary operator $ V $ is singular.

\begin{Thm} \label{thm5}
Let $\UU$ be a unitary operator such that $1$ is not an eigenvalue
of $\UU$. Then there exist
an inner function $\theta$ and an operator $\tS$
satisfying {\rm (i)--(iii)},
for which the restriction of $\tS$ to $\Kth$ is unitarily equivalent
to $\UU$ and $\phi_t(\tS)-\phi_t(S)\in\gS_p$ for all $p>1$.
\end{Thm}

For the case $p=2$, analogs of Theorems \ref{thm1}, \ref{thm5} were
obtained in \cite{AB2, AB3}. Here we also consider separately the results
about the class $\gS_2$, since in terms of our model their proofs are 
essentially simplified and one may expect that the conditions obtained 
are sharp. The new results of this paper are connected with more narrow classes, 
i.e., $\gS_p$ with $p<2$, and, in the first place, with $p=1$.

In \cite{AB,AB2} the construction of the operator $\tS$ was based on  
the triangulation studied by Ahern and Clark \cite{ac2}. Instead, we use 
Clark's construction \cite{clark}, which establishes an isometric
identification of the space $\Kth$ with a space $L^2(\mu)$ for a 
certain singular measure $\mu$ on the unit circle. 
This approach allows us to relate approximation properties of the difference 
$\phi_t(\tS)-\phi_t(S)$ with differential properties of $\mu$.

The situation becomes essentially different when we pass
to the ``natural'' unitary dilation of the perturbed semigroup $\tS$ with
properties (i)--(iii) (see \S2). 
Let $U$ be the bilateral shift on $L^2(\mathbb{T})$ and let $\tU$ be a dilation 
of $\tS$ with the Markovian property, i.e., such that $\tU^\ast$
coincides with $U^\ast$ 
on $H^2_-$. It turns out that then $\phi_t(\tU)-\phi_t(U)$ will never belong 
to $\gS_1$. However, for any unitary operator $\UU$ 
there exists a construction of $\tS$ with the properties {\rm (i)--(iii)} 
and a Markovian unitary dilation $\tU$ of $\tS$ such that
the restriction of $\tS$ to $\Kth$ is unitarily equivalent to $\UU$
and $\phi_t(\tU)-\phi_t(U) \in \gS_p$ for all $p>1$. 
Unitary dilations will be considered in detail elsewhere.

The paper is organized as follows.
It is shown in  \S2 that the conditions, imposed on the operator
$V$ in the main results, are necessary.
Then we consider at first the case where $\UU$ is a singular unitary operator 
of multiplicity 1 and we may apply Clark's construction. 
The model is determined by an 
inner function $\theta$, or by a singular measure $\mu$ on the
unit circle connected with $\theta$ by relation (\ref1) below, or by a measure
$\nu$ on the real line associated with $\mu$. 
In this case we have $rank\,(\tS-S)=1$. We obtain conditions on
$\mu$, $\nu$, under which the operators $\phi_t(\tS)-\phi_t(S)$
belong to $\gS_2$ and $\gS_1$ for all $t>0$.
We prove Theorem \ref{thm1} in the partial case of multiplicity 1
via properties of the measures $\mu$, $\nu$ from our construction,
and then from this special case we infer the general case of the theorem.
Theorems \ref{thm8} and \ref{thm5} are obtained as corollaries from 
Theorem \ref{thm1}.

The authors are grateful to R.V.~Romanov for the useful  discussions
and for the help with the proof of Theorem \ref{thm0}, 
and also to A.B.~Aleksandrov for many helpful suggestions which have improved
the exposition.

\bigskip

\section{Unitary groups, isometric semigroups, and their cogenerators}
\label{cogen}

In this section we recall some basic properties of 
generators and cogenerators of strongly continuous unitary 
groups and their analogs for isometric semigroups.
A detailed exposition of these topics may be found in 
\cite{io,sf,hf}.

If $(U_t)_{t\in\R}$ is a strongly continuous unitary group, 
then its {\it generator} $A$ is defined by 
$Ax = \lim_{t\to 0}\frac{U_t-I}{t} x$
on the set of those vectors $x$ for which the limit exists. 
Then $iA$ is a selfadjoint operator, and $U_t=\exp(tA)$.
The {\it cogenerator} is the unitary operator $B=(A+I)(A-I)^{-1}$.
A necessary and sufficient condition for a unitary operator to be 
a cogenerator of some unitary group is that its point spectrum does not 
contain the point 1. For the elements of the group we have the formula 
expressing them via the cogenerator:
\begin{equation}
\label{80a}
U_t=\phi_t(B)
\end{equation}
with $\phi_t$ defined by (\ref{80}). 
Moreover,
since 
$$
\int_0^{+\infty} e^{-t}\phi_t(z)\,dt=\frac{1-z}{2},
$$ 
we obtain
$$
\int_0^{+\infty} e^{-t}U_t\,dt=\int_0^{+\infty} e^{-t}\phi_t(B)\,dt=
\frac{I-B}{2};
$$
hence the cogenerator can be expressed via the elements of the group
with $t\geq 0$ as
$$
B=I-2 \int_0^{+\infty} e^{-t}U_t\,dt.
$$

Suppose now that we are given a strongly continuous semigroup $(V_t)_{t\ge 0}$
of isometric operators in a Hilbert space $H_0$. 
Then one may consider its unitary dilation, i.e., 
a unitary group $(U_t)_{t\in \mathbb{R}}$ acting on a wider space
$H$ such that $H_0$ is an invariant subspace for $U_t$ 
with $t>0$, and $V_t =U_t|_{H_0}$, $t>0$ (see \cite{sf}).
If $t\geq 0$, then $\phi_t\in H^\infty$. 
Therefore, $H_0$ is also an invariant subspace for 
the cogenerator $B$ of $(U_t)$, and it is natural to define the cogenerator of 
the isometric semigroup $(V_t)$ as the restriction of $B$
to $H_0$. Formula (\ref{80a}) holds true if instead of $B$ 
and $U_t$, $t>0$, we apply it to their restrictions to $H_0$.

For the classes of operators under consideration, the spectral type of 
an operator is defined as the class of all operators that 
are unitarily equivalent to the original one. 
Since unitary groups and isometric semigroups
are uniquely determined by their cogenerators and vice versa, 
it is natural to define the spectral type of a group (semigroup)
via the spectral type of its cogenerator.
For the class of unitary operators, 
the spectral type is determined by a scalar measure on the unit circle $\T$, 
under which the spectral measure of a unitary operator is absolutely 
continuous, and an integer-valued function on $\T$ counting the local 
multiplicity of the unitary operator at almost all points relative to 
the measure. According to the Wold--Kolmogorov decomposition of 
the cogenerator, every strongly continuous isometric semigroup splits 
into a direct sum of its unitary and pure isometric parts. The former is
a semigroup of unitary operators indexed by $\mathbb{R}_+$, and it
may be naturally extended to a unitary group with indices in $\mathbb{R}$.
The latter is a semigroup of completely non-unitary operators, i.e., operators
unitarily equivalent to a unilateral shift. Thus,
the spectral type of an isometric semigroup is determined by the spectral
type of the cogenerator of its unitary part and by the multiplicity of
the unilateral shifts.

The following (probably, known) statement shows that 
the multiplicity of the unilateral shift can not be changed by 
compact perturbations of the semigroup, while trace class perturbations 
preserve the absolutely continuous part of the cogenerator.
To make the exposition self-contained we present a short proof of this fact. 
The proof below was communicated to the authors
by R.V.~Romanov.

\begin{Thm} \label{thm0}
Let $(V_t)_{t\ge 0}$, $(\tilde V_t)_{t\ge 0}$  be strongly continuous
semigroups of isometric operators in a Hilbert space $H$.

$1)$ If the operator  $\tilde V_t- V_t$ is compact for all $t\ge 0$,
then the pure isometric parts of the semigroups $(V_t)_{t\ge 0}$ and 
$(\tilde V_t)_{t\ge 0}$ are unitarily equivalent, i.e.,
the multiplicities of the shifts coincide.

$2)$ If $\tilde V_t- V_t\in\gS_1$ for all $t\ge 0$, 
then the absolutely continuous parts of the semigroups
$(V_t)_{t\ge 0}$ and $(\tilde V_t)_{t\ge 0}$ are unitarily equivalent.
\end{Thm}

For the proof of the theorem we will need the following lemma.

\begin{Lem} \label{comp}
Let $T(t)$, $t\ge 0$, be a strongly continuous family of compact operators 
such that $\sup_{t\ge 0}\Vert T(t)\Vert<\infty$. Then the operator 
$\int_0^\infty e^{-t}T(t)dt$ is compact.
\end{Lem}

\beginpf  
For any $t$ consider the Schmidt expansion 
$T(t)=\sum_{j=1}^\infty s_{tj}(\cdot , x_{tj})y_{tj}$, 
where $x(_{tj})$, $(y_{tj})$ are orthonormal systems and
$s_{tj}\searrow 0$. Note that the functions $t\mapsto s_{tj}$ 
are continuous from below for all $j$ and, therefore, they are 
measurable. The functions $t\mapsto x_{tj}$ and $t\mapsto y_{tj}$.
are also measurable. The norm of the operator 
$\sum_{j\ge n} s_{tj}(h, x_{tj})y_{tj}$ equals $s_{tn}$, and hence 
$$
\Vert\int_0^\infty e^{-t} \sum_{j\ge n} s_{tj}(\cdot, x_{tj})y_{tj}\, dt\Vert
\leq \int_0^\infty e^{-t}s_{tn}\,dt.
$$
By the hypothesis, $\sup_{t\ge 0} s_{tn} \le \sup_{t\ge 0}\| T(t)\| <\infty$. 
For any $t$ we have $s_{tn}\to 0$, $n\to\infty$, whence,
by the Lebesgue theorem, we obtain $\int_0^\infty e^{-t}s_{tn}dt\to 0$, and so 
$\Vert\int_0^\infty e^{-t} \sum_{j\ge n} s_{tj}
(\cdot, x_{tj})y_{tj}\, dt\Vert \to 0$ as $n\to\infty$. 

It remains to show that the operator 
$\int_0^\infty e^{-t} \sum_{j<n} s_{tj}(\cdot, x_{tj})y_{tj}\, dt$
is compact for all $n$. To show this, note that 
the norms of the operators $(\cdot, x_{tj})y_{tj}$ in any class $\gS_p$ are equal 
to 1, whence $\gS_p$-norms of the operators  
$\int_0^\infty e^{-t} \sum_{j=1}^{n-1}s_{tj}(\cdot, x_{tj})y_{tj}\, dt$
do not exceed 
$$
\int_0^\infty e^{-t} \sum_{j=1}^{n-1} s_{tj}\,dt\leq (n-1)\cdot
\sup_{t\ge 0}\Vert T(t)\Vert \int_0^\infty e^{-t} dt =
(n-1)\cdot \sup_{t\ge 0}\Vert T(t)\Vert. 
$$
\qed 
\medskip
\\
\noindent {\bf Proof of Theorem \ref{thm0}.} \
We prove statement 1 of the theorem. Applying Lemma \ref{comp} to 
$T(t)=\tilde V_t- V_t$, we obtain that the operator 
$\int_0^\infty e^{-t}(\tilde V_t- V_t)dt$ is compact.
This means that the difference of the cogenerators 
of the semigroups $( V_t)$ and  $(\tilde V_t)$ is compact. 
Hence, the cogenerators have equal Fredholm indices, which are exactly 
the multiplicities of the shift.

Now we prove statement 2. Denote by $Q$ the subspace where
the cogenerator of the semigroup $(\tilde V_t)$ acts as an absolutely continuous 
unitary operator. Denote by $Z$ the natural embedding of the
$Q$ into $H$. It follows from the assumption   
$\tilde V_t- V_t\in\gS_1$ that 
$Z (\tilde V_t|_Q)-\gamma_t Z\in\gS_1$, where $(\gamma_t)_{t\in\R}$ is 
the unitary dilation of the semigroup $( V_t)$. 
Then, by the classical scattering theory for a pair
of unitary operators, a strong limit 
of the isometries $\gamma_t Z(\tilde V_{t}|_Q)^{-1}$ exists as  $t\to +\infty$,
which defines an isometric wave operator $W$ (see \cite[Theorem 6.5.5]{yaf}). 
By the construction, the range of the operator $W$ is contained in $H$ 
and reduces the group $(\gamma_t)$. Hence, the operator 
$W$ realizes a unitary equivalence
between the restriction of the semigroup $(\tilde V_t)$ to $Q$ 
and some unitary part of the semigroup $(V_t)$.

We have shown that the absolutely continuous unitary part
of the semigroup $(\tilde V_t)$ is unitary equivalent to 
some part of the semigroup $( V_t)$.
Analogously, the absolutely continuous unitary part
of the semigroup  $( V_t)$ is unitarily equivalent to 
some part of the semigroup $(\tilde V_t)$. 
Then, by the spectral theorem, the absolutely continuous unitary parts 
of the semigroups $( V_t)$ and $(\tilde V_t)$ are unitarily equivalent. \qed
\medskip

We have an immediate consequence for the perturbations of the semigroup 
of shifts.

\begin{Cor} \label{neccond}
Let $(\tau_t)_{t\ge 0}$ be the semigroup of shifts and let
$(\tilde\tau_t)_{t\ge 0}$ be a strongly continuous semigroup 
of isometric operators in $L^2(\R_+)$. 

$1)$ If the operator $\tilde\tau_t-\tau_t$ is compact for all 
$t\ge 0$, then there exists a unitary group $(\omega_t)_{t\in\R}$ such that  
the semigroup $(\tilde\tau_t)$ is unitarily equivalent to the direct sum 
$(\tau_t)\oplus(\omega_t)$, $t\ge 0$. 

$2)$ If, additionally, $\tilde\tau_t-\tau_t\in\gS_1$
for $t\ge 0$, then the  spectral measure $($of each element 
or, equivalently, of the cogenerator$)$ of the group $(\omega_t)$ is singular.
\end{Cor}

For the model on the unit circle Corollary \ref{neccond} 
means that under the assumption $\phi_t(\tS)-\phi_t(S)\in\gS_1$, $t>0$, 
the multiplicity of the shift for the operator $\tS$
is 1, and the spectral measure of its unitary part is necessarily singular 
with respect to the Lebesgue measure.

\bigskip

\section{The spaces $\Kth$ and the model construction}\label{spaces}

In this section we introduce a special model of a perturbation 
satisfying {\rm (i)--(iii)} starting from a singular measure $\mu$ 
on the unit circle $\T$ and an inner function $\theta$ 
associated with $\mu$.

Let $\mu$ be a measure on $\T$ singular relative to 
the Lebesgue measure $m$, and let $\mu(\{1\}) = 0$. 
Define the function $\theta$ by 
\begin{equation}
\label{1}
\frac{1+\theta(z)}{1-\theta(z)}=\int_\T \frac{1+\bar\xi z}{1-\bar\xi z}
\, d\mu(\xi), \qquad z\in\mathbb{D}.
\end{equation}
It is well known that $\theta$ is an inner function.
The measure $\mu$ is concentrated on the set where angular boundary limits
of $\theta$ exist and are equal to 1. The measure $\mu$ is said to be a Clark 
measure of the function $\theta$. 

For $u\in L^2(\mu)$ put
\begin{equation}
\label{46}
(\Om u)(z)=(1-\theta(z))\int_\T\frac{u(\xi)d\mu(\xi)}{1-\bar\xi z}.
\end{equation}
Clark  \cite{clark} proved that $\Om$
is a unitary operator from $L^2(\mu)$
to $\Kth$. Moreover, angular boundary values of the function $\Om u$ exist
and coincide with $u$\; $\mu$-almost everywhere \cite{polt}. Analogs of
these results for the spaces of vector-valued functions may be found
in \cite{kp}. 

In this section $\UU$ is the operator of multiplication by the independent 
variable $\xi$ on $L^2(\mu)$. Note that $1$ is not an eigenvalue for
$\UU$, since $\mu(\{1\}) = 0$. We find a formula for the unitary
operator $\Om\UU\Om^\ast$, which is a unitarily equivalent transplantation 
of $\UU$ to $\Kth$. For $h=\Om u$, $u\in L^2(\mu)$, we have
$$
\begin{aligned}
  (\Om\UU\Om^\ast h)(z)-zh(z)
  & = (\Om\UU u)(z)-z(\Om u)(z) \\
  & = (1-\theta(z))\int_\T \frac{(\xi-z)u(\xi)d\mu(\xi)}{1-\bar\xi z} \\
  & = (1-\theta(z))\int_\T \xi u(\xi)d\mu(\xi).
\end{aligned}
$$
Since 
$\int_\T \xi u(\xi)d\mu(\xi)=(u, \bar\xi)_{L^2(\mu)}=(h, \Om\bar\xi)_{\Kth}$, 
we obtain
$$
\Om\UU\Om^\ast h=zh+(h, g)(1-\theta), \quad h\in\Kth,
$$
where $g=\Om\bar\xi\in\Kth$ (it is easy to check the formula for $g$,
$g(z)=\frac{\theta(z)-\theta(0)}{z(1-\theta(0))}$, which will not be
used in this paper).

Recall that $S$ denotes the shift operator on $H^2$, and define 
the operator $\tS$ on $H^2$ by 
\begin{equation}
\label{2}
\tS=S+(\cdot, g)(1-\theta).
\end{equation}
We have shown that $\Kth$ is an invariant subspace for $\tS$,
and that the restriction of $\tS$ to $\Kth$ is unitarily equivalent to 
$\UU$. Clearly, $\tS$ coincides with $S$ on $\Kth$. We have obtained
the following proposition.

\begin{Prop} For the operator $\tS$ defined by $(\ref{2})$
properties {\rm (i)--(iii)} are fulfilled; the restriction of $\tS$
to $\Kth$ is unitarily equivalent to the operator $\UU$ of multiplication
by the independent variable on the space $L^2(\mu)$.
Thus, $\tS$ is unitarily equivalent to $S\oplus\UU$.
\end{Prop}

The operator $\tS$ differs from the multiplication by $z$ by a rank-one 
operator:  $\tS-S=(\cdot, g)(1-\theta)$, for the norm of which we have
\begin{equation}
\label{47}
  \|\tS-S\|=\|g\|_{\Kth}\cdot \|1-\theta\|_{H^2}
   = \|\bar\xi\|_{L^2(\mu)}\cdot \|1-\theta\|_{H^2}<2\sqrt{\mu(\T)}.
\end{equation}

We shall work mainly with measures $\mu$ satisfying an additional condition
\begin{equation}
\label{3}
\int_\T \frac{d\mu(\xi)}{|1-\xi|^2}<\infty,
\end{equation}
This implies, in particular, that $\mu(\{1\}) =0$. 
Condition (\ref{3}) is well known in the theory of inner functions
and model subspaces. It is equivalent to any of the following 
(see, e.g., \cite[Chapter VI]{sar}):  

(i) the function $\theta$ defined by (\ref{1}) 
has a finite angular derivative at the point 1;

(ii) each function from $K_\theta$ has a finite
nontangential limit at the point 1.

Moreover, the function $\frac{1-\overline{\theta(1)}\theta}{1-z}$
belongs to $\Kth$ and is the reproducing kernel at the point 1,
$$
\Big\|\frac{1-\overline{\theta(1)}\theta}{1-z}\Big\|^2_{\Kth}=
\Big\|\frac{1-\overline{\theta(1)}}{1-z}\Big\|^2_{L^2(\mu)}
< 4 \int_\T \frac{d\mu(\xi)}{|1-\xi|^2}.
$$
\bigskip

\section{Equivalent models}\label{eqmodels}

Our original construction was a model of shifts in
the space $L^2(\R_+)$. We shall also work in other models, which
are unitarily equivalent the original one.

An equivalent model construction on the real line may be 
obtained from the model of shifts by means of 
the Fourier transform. It sends the space $L^2(\R_+)$ to
the Hardy space $H^2(\C_+)$ in the upper half-plane.
The isometric semigroup $(\tau_t)_{t\ge 0}$ of shifts becomes the semigroup
of operators of multiplication by the functions $\exp(itz)$.

For the model on the unit circle $\T$ we work in the Hardy space $H^2$,
which is a subspace of $L^2$. In the notations $L^2, H^2$ we omit
the measure, then this means  the Lebesgue measure
$m$ on $\T$ normalized so that $m(\T)=1$.
The cogenerator of the unperturbed isometric semigroup is the operator $S$
of multiplication by the independent variable $z$. Then the semigroup
consists of operators of the form $\phi_t(S)$ with $\phi$ defined by
(\ref{80}). In other words, our isometric semigroup is the group of
operators of multiplication by the inner functions $\phi_t$, $t\geq 0$.

Now we list some formulas establishing a unitary equivalence of the models
of multiplications on the unit circle $\T$ and on the real line $\R$.
For the variable $z$ on $\T$ we write $x=i\frac{1+z}{1-z}\in\R$.
Given a measure $\mu$ on $\T$, define the measure $\nu$ on $\R$ by 
\begin{equation}
\label{meas-trans}
d\mu (z) = \frac{d\nu(x)}{\pi(1+x^2)}.
\end{equation} 
Condition (\ref{3}) is equivalent to $\nu(\R)<\infty$.
The mapping 
\begin{equation}
\label{uv}
u\mapsto v, \quad v(x)=\frac{1}{\sqrt{\pi}(x+i)}
\cdot u\Big(\frac{x-i}{x+i}\Big),
\end{equation}
is a unitary operator from $L^2(\mu)$ onto $L^2(\nu)$, and also from
$L^2 = L^2(\mathbb{T})$ to $L^2(\R)$; the latter maps the Hardy class $H^2$
onto $H^2(\C_+)$. The function $u\in L^2(\mu)$ is expressed via
$v\in L^2(\nu)$ as
$$
u(z)=\frac{2i\sqrt\pi}{1-z}\cdot v\bigg(i\frac{1+z}{1-z}\bigg).
$$

\bigskip

\section{Functions of operators $S$ and $\tS$}\label{functions}

We deal with the construction,
where $\tS$ is a rank-one perturbation of $S$ defined by (\ref2).
Take a function $\phi\in H^\infty$ and suppose that values of $\phi$
are defined $\mu$-almost everywhere.

On the subspace $\theta H^2$, the operator $\tS$ 
coincides with $S$, hence $\phi(\tS)$
coincides with $\phi(S)$ there. Thus, we need to study only the restriction
of $\phi(\tS)-\phi(S)$ to $\Kth$, which makes natural to consider
the operator
\begin{equation}
\label{operx}
X: L^2(\mu)\to H^2, \quad X=(\phi(\tS)-\phi(S))\Om,
\end{equation}
with $\Omega$ defined by (\ref{46}). Take $u\in L^2(\mu)$. Then
\begin{equation}
\label{xu}
\begin{aligned}
(Xu)(z) & = (1-\theta(z))\int_\T\frac{\phi(\xi)u(\xi)}{1-\bar\xi z}d\mu(\xi) 
        - \phi(z)(1-\theta(z))\int_\T\frac{u(\xi)}{1-\bar\xi z}d\mu(\xi) \\
        & = (1-\theta(z))\int_\T\frac{\phi(\xi)-\phi(z)}{\xi-z}\xi
u(\xi)d\mu(\xi), \qquad u\in L^2(\mu).
\end{aligned}
\end{equation}
Since $\Omega$ is a unitary operator from $L^2(\mu)$ onto $\Kth$,
the operator $\phi(\tS)-\phi(S)$ belongs to a class $\gS_p$ if and only
if $X\in\gS_p$, and their $\gS_p$-norms coincide.

Equality (\ref{uv}) establishes unitary correspondences 
between $L^2(\mu)$ and $L^2(\nu)$, and between $H^2$ and $H^2(\C_+)$. 
Using these identifications, we construct the operator $Y: L^2(\nu)\to H^2(\C_+)$ as
a unitary transplantation of $X$. We have
\begin{equation}
\label{opery}
\begin{aligned}
(Yv)(x) & =\frac{1}{\sqrt{\pi}(x+i)}(Xu)\left(\frac{x-i}{x+i}\right) \\
& =\frac{1}{\sqrt{\pi}(x+i)}\cdot (1-\theta(z))
\int_\T\frac{\phi(\xi)-\phi(z)}{\xi-z}\xi u(\xi)d\mu(\xi) \\
& =(1-\Theta(x))\cdot\frac{1}{2\pi i}\int_\R\frac{\psi(\zeta)-\psi(x)}{\zeta-x}
v(\zeta)d\nu(\zeta), \qquad v\in L^2(\nu),
\end{aligned}
\end{equation}
where 
$$
\psi(x)=\phi\left(\frac{x-i}{x+i}\right), \qquad
\Theta(x)=\theta\left(\frac{x-i}{x+i}\right), \qquad x\in \mathbb{C}_+.
$$
Note that if $\theta$ is defined by (\ref{1}), then for $\Theta$ we have
the formula
$$
\frac{1+\Theta(x)}{1-\Theta(x)}=\frac{1}{\pi i}\int
\left(\frac{1}{\zeta-x}-\frac{\zeta}{1+\zeta^2}\right)d\nu(\zeta), 
\qquad x\in \mathbb{C}_+.
$$

The following proposition  follows directly from the construction of $Y$.

\begin{Prop} \label{yy}
The operators $\phi(\tS)-\phi(S)$ and $Y$ belong or not to ideals 
$\gS_p$ simultaneously, and $\|\phi(\tS)-\phi(S)\|_{\gS_p}=\|Y\|_{\gS_p}$.
\end{Prop}
\bigskip

\section{Estimates for the Hilbert--Schmidt class}\label{hs}

Define the integral operator $K: L^2(\nu)\to H^2(\C_+)$ by
\begin{equation}
\label{62}
(Kv)(x)=\frac{1}{2\pi i}\int\frac{\psi(\zeta)-\psi(x)}{\zeta-x}
v(\zeta)d\nu(\zeta), \qquad v\in L^2(\nu).
\end{equation}
If $K\in\gS_p$, then $Y\in\gS_p$ as well, and we obviously have 
\begin{equation}
\label{444}
\|Y\|_{\gS_p}\leq 2\cdot\|K\|_{\gS_p}.
\end{equation}

For the norm of $K$ in $\gS_2$ we obtain
\begin{equation}
\label{41}
\|K\|_{\gS_2}^2 =\frac{1}{4\pi^2} 
\iint_{\R\times\R}\Big|\frac{\psi(\zeta)-\psi(x)}{\zeta-x}\Big|^2
dx\,d\nu(\zeta).
\end{equation}

Consider the semigroups of the form $(\phi_t(S))$, $(\phi_t(\tS))$,
where $\phi_t(z)=\exp(t\frac{z+1}{z-1})$, $t>0$. The corresponding
functions $\psi_t$ have the form
$$
\psi_t(x)=\phi_t\left(\frac{x-i}{x+i}\right)=e^{itx}.
$$
Now we state a condition on $\nu$ ensuring that the operator $K$ belongs
to the Hilbert--Schmidt class for all $t$ at once. To this
end, we shall obtain a precise formula for the integral in (\ref{41})
with $\psi=\psi_t$; it will give us a sufficient 
condition for $\phi_t(\tS)-\phi_t(S)\in\gS_2$ for all $t$.
At the end of the section we present
an example of the construction, for which $\phi_t(\tS)-\phi_t(S)$ is not
in the Hilbert-Schmidt class for any $t\ne 0$.

\begin{Prop} \label{gilb_class}
Assume that $(\ref{3})$ is fulfilled. Then for any $t>0$ we have
\begin{equation}
\label{4}
\iint_{\R\times\R}\Big|\frac{\psi_t(\zeta)-\psi_t(x)}{\zeta-x}\Big|^2dx\,
d\nu(\zeta)=2\pi t\cdot\nu(\R)=8\pi^2 t
\int_\T \frac{d\mu(\xi)}{|1-\xi|^2}<\infty.
\end{equation}
\end{Prop}

\beginpf 
For the integral on the left-hand side with respect to the Lebesgue measure, 
we obtain
$$
\int_\mathbb{R} \frac{|e^{it\zeta}-e^{itx}|^2}{(\zeta-x)^2}\, dx
=4 \int_\mathbb{R}\frac{\sin^2 \frac{t}{2}(\zeta-x)}{(\zeta-x)^2} \,dx
=2\pi t,
$$
where we used the well-known formula
$\int_\mathbb{R}\frac{\sin^2\alpha s}{s^2}ds=\pi |\alpha|$.
The right identity in (\ref{4}) follows from the relation
$\nu(\R)=4\pi\int\frac{d\mu(\xi)}{|1-\xi|^2}$.
\qed

\begin{Cor}
If $\mu$ satisfies $(\ref{3})$, then $\phi_t(\tS)-\phi_t(S) \in \gS_2$
for all $t>0$;
\begin{equation}
\label{99}
\|\phi_t(\tS)-\phi_t(S)\|_{\gS_2}\leq 2\sqrt{2t}
\left(\int_\T\frac{d\mu(\xi)}{|1-\xi|^2}\right)^{1/2}.
\end{equation}
\end{Cor}

\beginpf This directly follows from Proposition \ref{yy} and relations
(\ref{444}), (\ref{41}), and (\ref{4}).
\qed

\medskip

In \S7 we find a certain ``smallness'' condition for
a measure $\mu$ at the point 1, see (\ref{tr}), which is sufficient 
for the inclusion of the operator $\phi_t(\tS)-\phi_t(S)$ 
into the trace class $\gS_1$. 

From the representation (\ref{opery}) of $Y$ in the form of an integral
operator it follows that $Y\in\gS_2$ if and only if
$$
\iint_{\R\times\R}  |1-\Theta(z)|^2
\Big|\frac{\psi(\zeta)-\psi(x)}{\zeta-x}\Big|^2 dx d\nu(\zeta)<\infty.
$$
We expect that the condition $\nu(\R)<\infty$ equivalent to (\ref{3}),
is also necessary for the inclusions $\phi_t(\tS)-\phi_t(S) \in \gS_2$,
that is, the factor $1-\Theta(z)$
does not change essentially the convergence of the integral in (\ref{41}).

Now we present an example of a measure $\nu$ such that the corresponding
operator $\phi_t(\tS)-\phi_t(S)$ is not in the Hilbert--Schmidt class
for any $t\ne 0$. Recall that $\tS-S$ is a rank-one operator whose
norm can be made arbitrarily small.

\begin{Ex}
{\rm Let $\nu = \sum_{n\in \mathbb{Z}} \delta_n$ be the sum of unit point
masses at all integers. Thus, the total mass of $\nu$ is infinite, and hence
the corresponding measure $\mu$ on $\T$ does not satisfy condition (\ref{3}).
For $\psi_t(x)=e^{itx}$ in the last formula, we obtain the integral
$$
\int_\R |1-\Th(s)|^2\, \bigg(\sum_{n\in\mathbb{Z}}
\frac{\sin^2 \frac{t}{2}(s-n)}{(s-n)^2} \bigg).
$$
Obviously, for any $t\ne 0$, there exist positive constants $\delta(t)<1/2$
and $C(t)$ such that the sum in brackets in the above formula 
is bounded from below by $C(t)$ for $s\in (k+\delta(t), k+2\delta(t))$, 
$k\in\mathbb{Z}$. On the other hand, an elementary estimate shows that
for $s\in (k+\delta(t), k+2\delta(t))$, $k\in\mathbb{Z}$,
$$
\bigg|\frac{1+\Th(s)}{1-\Th(s)}\bigg|  =
\bigg|\sum_{n\in\mathbb{Z}}  \frac{1+ s n}{(s-n)(n^2+1)} \bigg| \le C_1(t),
$$ 
and therefore $|1-\Theta(s)|\ge C_2(t)>0$. Combining
these estimates, we conclude that the integral diverges, i.e., $Y$ does
not belong to the Hilbert--Schmidt class. } 
\end{Ex}
\bigskip

\section{Estimates for the trace class}\label{trace}

In this section, instead of the Hilbert--Schmidt class $\gS_2$, we consider
the classes $\gS_p$ with $p<2$, and, in the first place, 
the trace class $\gS_1$.
According to (\ref{444}), to show that $\phi_t(\tS)-\phi_t(S) \in \gS_p$,
we need to prove that $K\in\gS_p$, where $K$ is defined by (\ref{62})
with $\psi=\psi_t$, $\psi_t(x)=e^{itx}$. We reduce this problem
to a question about embeddings of the Paley--Wiener space $\PW_t$ of all
entire functions of exponential type at most $t$ that are square-summable
on $\mathbb{R}$. Embedding operators for spaces 
of analytic functions and their inclusion to ideals $\gS_p$ were studied 
in detail by O.G.~Parfenov \cite{parf2, parf}. In \cite{bar08} 
generalizations of some of Parfenov's results to embeddings of 
coinvariant subspaces $K_\theta$ are obtained.

Consider the operator adjoint to the operator 
$K$ defined by (\ref{62}) with $\psi=\psi_t$, $\psi_t(x)=e^{itx}$:
$$
\begin{aligned}
(K^\ast f)(\zeta) & =-\frac{1}{2\pi i}
\int\frac{e^{-it\zeta}-e^{-itx}}{\zeta-x}f(x)dx \\
& =e^{-\frac{it\zeta}2}\int e^{-\frac{itx}2}f(x)
\left(\frac{1}{\pi}\cdot\frac{\sin\frac{t}{2}(x-\zeta)}{x-\zeta}\right)dx,
\quad f\in L^2(\R).
\end{aligned}
$$
It is well known that the function
$\frac{1}{\pi}\cdot\frac{\sin\frac{t}{2}(x-\zeta)}{x-\zeta}$ is
the reproducing kernel at $\zeta$ for the space $\PW_{t/2}$, which
means that this function belongs to $\PW_{t/2}$, and the inner
product of any $f \in \PW_{t/2}$ with the reproducing kernel
equals $f(\zeta)$. Therefore,
$K^\ast f=0$ if the function $\tilde f$,
$\tilde f(x)=e^{-\frac{itx}2}f(x)$, is orthogonal to $\PW_{t/2}$,
and if $\tilde f\in \PW_{t/2}$, then
$$
(K^\ast f)(\zeta)=e^{-\frac{it\zeta}2}\cdot \tilde f(\zeta)
=e^{-it\zeta}f(\zeta).
$$
Thus, $K^\ast$ belongs to a class $\gS_p$ simultaneously with
the operator $E_{\nu, t/2}$ which embeds the space $\PW_{t/2}$
into $L^2(\nu)$. 

The following criterion is due to Parfenov \cite{parf2, parf}.

\begin{Thm}
Let $\Delta_n = [n,n+1)$. Then $E_{\nu, t}\in\gS_p$, $p>0$, if and only if
\begin{equation}
\label{criter}
\sum\limits_{n\in \Z} \big(\nu(\Delta_n)\big)^{p/2} <\infty.
\end{equation} 
Moreover,
\begin{equation}
\label{criter1}
\|E_{\nu, t}\|_{\gS_p}^p\le C(p)\, t^{p/2}
\sum\limits_{n\in \Z} \big(\nu(\Delta_n)\big)^{p/2}.
\end{equation} 
\end{Thm}

Going back to the measure $\mu$, condition (\ref{criter}) reads as 
$$
\sum\limits_{n\in \Z} \bigg(\int_{\gamma_n} \frac{d\mu(\xi)}{|1-\xi|^2} 
\bigg)^{p/2} <\infty,
$$
where the arcs $\gamma_n$ are defined for $n>0$ by 
$\gamma_n = \{e^{i \phi}: \pi/(n+1) \le \phi\le \pi/n\}$ and
symmetrically for $n<0$.

We obtain a condition, which is analogous to (\ref{3}) 
and is sufficient for (\ref{criter}) and, therefore, sufficient
for the inclusion of $\phi_t(\tS)-\phi_t(S)$ into the class
$\gS_p$.

\begin{Prop} 
\label{tr_class}
Let $0<p<2$. If the measure $\mu$ satisfies
\begin{equation}
\label{tr}
\int_{\T} \frac{d\mu(\xi)}{|1-\xi|^q} <\infty
\end{equation} 
for some $q> 1+ 2/p$, then $K\in \gS_p$, and, consequently,
$\phi_t(\tS)-\phi_t(S) \in \gS_p$.

In particular, $(\ref{tr})$ with $q>3$ yields $\phi_t(\tS)-\phi_t(S) \in \gS_1$
and
\begin{equation}
\label{996}
\|\phi_t(\tS)-\phi_t(S)\|_{\gS_1}\le
M_q\cdot t^{1/2}\cdot\bigg(\int_\T \frac{d\mu(\xi)}{|1-\xi|\,^q}\bigg)^{1/2},
\end{equation}
where $M_q$ is a constant depending only on $q$.
\end{Prop} 

\beginpf 
Rewrite condition (\ref{tr}) in terms of the measure $\nu$ defined
by (\ref{meas-trans}):
\begin{equation}
\label{tr-up}
\int_\R (x^2+1)^{r/2}d\nu(x)=2^{-q}\pi\int_{\T}\frac{d\mu(\xi)}{|1-\xi|^q}
<\infty,
\end{equation}
where $r=q-2> (2-p)/p$. Then, by the H\"older inequality with exponents
$2/p$ and $2/(2-p)$, 
$$
\begin{aligned}
\bigg(\sum\limits_{n\in \Z} \big(\nu(\Delta_n)\big)^{p/2}\bigg)^{2/p}  & \le
\sum\limits_{n\in \Z} (|n|+1)^{r}\, \nu(\Delta_n)  \cdot
\bigg(\sum\limits_{n\in \Z} (|n|+1)^{-pr/(2-p)}\bigg)^{(2-p)/p} \\
& = {\rm const} \cdot \sum\limits_{n\in \Z} (|n|+1)^{r}\, \nu(\Delta_n) \le
{\rm const} \cdot \int_\R (|t|+1)^{r} d\nu (t).
\end{aligned}
$$
Now the statement follows from Parfenov's theorem and inequality
(\ref{criter1}).
\qed

\begin{Ex}
{\rm The exponent 3 is sharp, and for $q=3$ the second statement of
Proposition \ref{tr_class} fails. Indeed, choose $\nu$ so that
$\nu (\Delta_n)=\big((|n|+1)\log (|n|+2)\big)^{-2}$. Then (\ref{tr-up})
is fulfilled, but (\ref{criter}) is not, and hence $K$ is not in $\gS_1$.}
\end{Ex}
\bigskip

\section{Operators of multiplicity $>1$}\label{multip}

Now we use the construction developed above for the case of multiplicity 1,
for a generalization to the operators with arbitrary multiplicity. 

Let $\{\mu_n\}$ be a family of singular measures on the unit circle, 
where $n$ runs over the set $n=1, 2, \dots, N$ for some positive integer
$N$ or over $\mathbb{N}$. We assume that for some $q>3$
\begin{equation}
\label{93}
\sum_n\left(\int_\T\frac{d\mu_n(\xi)}{|1-\xi|^q}\right)^{1/2}<\infty.
\end{equation}

For each $n$, we construct objects as in \S3: 
the function $\theta_n$ is determined by formula (\ref{1}) with $\mu_n$ 
in place of $\mu$, the operator $\Om_n: L^2(\mu_n)\to K_{\theta_n}$
acts by formula (\ref{46}) with $g_n=\Om_n\bar\xi$, $\UU_n$ is the operator
of multiplication by $z$ on $L^2(\mu_n)$. Set 
$$
\hat\theta_n=\prod_{k=1}^{n-1}\theta_k.
$$
The operator $\Om: \sum\oplus L^2(\mu_n)\to H^2$ will be defined here by
$$
\Om\Big(\sum\oplus u_n\Big)=\sum\hat\theta_n\Om_n u_n. 
$$
Since condition (\ref{93}) is fulfilled for some $q>3$, it is also true for
$q=2$. Hence,
\begin{equation}
\label{812}
\sum_n \Big\|\frac{1-\overline{\theta_n(1)}\theta_n}{1-z}\Big\|_{L^2}
<2\sum_n\left(\int_\T \frac{d\mu_n(\xi)}{|1-\xi|^2}\right)^{1/2}<\infty.
\end{equation}
Therefore, the series 
$$
\sum_n\overline{\hat\theta_n(1)}\hat\theta_n
\frac{(1-\overline{\theta_n(1)}\theta_n)}{1-z}
$$
converges in the norm of the space $L^2$. The partial sums have
the form $\frac{1-\overline{\hat\theta_n(1)}\hat\theta_n}{1-z}$. Hence the 
limit function can be written as $\frac{1-\theta}{1-z}$ for an inner function 
$\theta$ that coincides with the product of the functions $\theta_n$ up to 
a unimodular multiplicative constant. It is easily seen that $\Om$ 
isometrically maps $\sum\oplus L^2(\mu_n)$ onto $\Kth$.

Define the operator $\tS$ on $H^2$ by
\begin{equation}
\label{70}
\tS=S+\sum_n(\cdot, \hat\theta_n g_n)\,\hat\theta_n\cdot(1-\theta_n),
\end{equation}
which is a generalization of (\ref{2}).
Then the operator $\tS$ is diagonal with respect to the decomposition 
$H^2=\sum\oplus \hat\theta_n K_{\theta_n}\oplus\theta H^2$. 
Indeed, $\tS$ coincides with $S$ on $\theta H^2$, while on each subspace
$\hat\theta_n K_{\theta_n}$, $\tS$ is a unitary transplantation of
multiplication by $z$ on $L^2(\mu_n)$. 
By a direct computation, we obtain  
$\tS(\hat\theta_n f)=\hat\theta_n \Omega_n (z u)\in \hat\theta_n K_{\theta_n}$
for $f=\Omega_n u \in K_{\theta_n}$.

Thus, $\tS$ satisfies properties (i)--(iii).
Similarly to (\ref{47}), we obtain the estimate for the trace class norm of
$\tS-S$:
\begin{equation}
\label{212}
\|\tS-S\|_{\gS_1}<2\sum_n\sqrt{\mu_n(\T)}.
\end{equation}
To estimate the trace class norm of $\phi_t(\tS)-\phi_t(S)$, we
use inequality (\ref{996}), which gives us
\begin{equation}
\label{81}
\|\phi_t(\tS)-\phi_t(S)\|_{\gS_1}\leq M_q\cdot t^{1/2}\cdot
\sum_n \left(\int_\T \frac{d\mu_n(\xi)}{|1-\xi|\,^q}\right)^{1/2}.
\end{equation}

\medskip

\noindent {\bf Proof of Theorem \ref{thm1}.} \
Take an arbitrary unitary operator $\UU$ whose spectral measure is singular
relative to the Lebesgue measure and has no point mass at $1$. 
Then there exist singular unitary operators $\UU_n$ of multiplicity 1 such
that $\UU=\oplus\sum\UU_n$. Construct measures $\mu_n$ so that each operator 
$\UU_n$ is unitarily equivalent to multiplication by the independent variable 
$\xi$ on $L^2(\mu_n)$.
We may think the measures $\mu_n$ satisfy (\ref{93}) with $q>3$ and 
$\sum_n\sqrt{\mu_n(\T)} < \eps/2$; otherwise, these properties may be 
fulfilled after multiplying the measures $\mu_n$ by appropriate positive weights.
Consider the operator $\tS$ defined by (\ref{70}).
Then (\ref{212}) yields $\|\tS-S\|_{\gS_1}<\eps$, and for any $t>0$ the operator 
$\phi_t(\tS)-\phi_t(S)$ is of trace class with norm estimated by (\ref{81}).
\qed

\begin{Rem}
{\rm For the Hilbert--Schmidt norm, from (\ref{99}) we obtain the estimate}
$$
\|\phi_t(\tS)-\phi_t(S)\|_{\gS_2}\leq 2\sqrt{2t}
\left(\sum_n\int_\T\frac{d\mu_n(\xi)}{|1-\xi|^2}\right)^{1/2}.
$$
\end{Rem}

\noindent {\bf Proof of Theorem \ref{thm8}.}
Statement 1 is a particular case of Theorem \ref{thm0}. Statement  2
follows immediately from  Theorem \ref{thm1} 
when we pass to the model connected with the group 
of shifts on the line. Indeed, in the construction in Theorem 
\ref{thm1} the differences of the elements of the semigroups 
on $L^2(\R_+)$ belong to the trace class, as required. \qed

\medskip
                    
\noindent {\bf Proof of Theorem \ref{thm5}.} 
Let $V$ be an arbitrary unitary operator such that 1 is not its eigenvalue.
Then $V=(A-iI)(A+iI)^{-1}$ for some selfadjoint operator $A$, and
$\phi_t(V)=\exp(itA)$. Applying a variant of the Weil--von Neumann
theorem which is due to Kuroda  (see \cite{kato} and \cite[Theorem
6.2.5]{yaf}) to some cross-normed ideal which is contained in all 
ideals $\gS_p$ with $p>1$, 
we may construct a close selfadjoint operator with pure point spectrum  
for a given selfadjoint operator 
so that the difference between the perturbed operator
and the original one belongs to  $\gS_p$ for all
$p>1$, and the norms may be taken to be arbitrarily small.

We represent the operator $A$  as the direct sum 
of bounded operators  $A_n$ and, applying the Kuroda theorem,
construct operators $A_n'$, so that the norms $\|A_n-A_n'\|_{\gS_p}$ are
small. It follows from the results by Davies \cite{davis} that 
$\exp(itA_n)-\exp(itA'_n)$ belongs to all classes
$\gS_p$, $p>1$, for $0\le t\le 1$. Thus, we
may construct a selfadjoint operator $A'$ 
such that $\exp(itA)-\exp(itA')$ is in all classes
$\gS_p$, $p>1$, for $0\le t\le 1$. Hence this holds for any 
$t>0$, since the inclusion $\exp(itA)-\exp(itA') \in \gS_p$
implies the inclusion $\exp(2itA)-\exp(2itA') \in \gS_p$. 
Now the statement follows from Theorem \ref{thm1} 
applied to the unitary operator $V'=(A'-iI)(A'+iI)^{-1}$.
\qed

\renewcommand{\refname}{References}

\enddocument